\documentclass[a4paper, 12pt]{article}

\usepackage[koi8-r,cp1251]{inputenc}
\usepackage{amsfonts,amssymb,amsmath, hyperref}
\usepackage[final]{epsfig}
\usepackage{graphicx}

\begin{document}

\begin{center}
\textbf{\large{ On the essential spectrum of $\lambda$-Toeplitz operators over compact Abelian groups}
}
\end{center}

\

\begin{center}
{A. R. Mirotin}
\end{center}

\

\begin{center}
amirotin@yandex.ru
\end{center}

\

\textsc{
Abstract.} 
In the recent paper by Mark C. Ho (2014) the notion of  a $\lambda$-Toeplitz operator on the Hardy space $H^2(\mathbb{T})$ over the one-dimensional torus $\mathbb{T}$ was introduced and it was shown  (under the supplementary condition) that for  $\lambda\in \mathbb{T}$ the essential spectrum of such an operator  is invariant with respect to the rotation $z\mapsto \lambda z$; if in addition $\lambda$ is not of finite order the essential spectrum is circular.  In this paper, we generalize these results to the case when  $\mathbb{T}$ is replaced by an arbitrary  compact Abelian group whose  dual is totally ordered.

\


\

MSC [2008] 47B35,  47A17

\section{Introduction}
\label{1}

In the paper \cite{Ho}   the notion of a $\lambda$-Toeplitz operator on the Hardy space $H^2(\mathbb{T})$ over the one-dimensional torus $\mathbb{T}$ was introduced and it was shown (under the supplementary  condition that  some modification of the symbol  belongs to $L^\infty(\mathbb{T})$) that for  $\lambda\in \mathbb{T}$ the essential spectrum of such an operator  is invariant with respect to the rotation $z\mapsto \lambda z$; if in addition $\lambda$ is not of finite order the essential spectrum is circular.  In this paper, we generalize these results to the case when the group $\mathbb{T}$ is replaced by an arbitrary  compact Abelian group whose  character group is totally ordered. We  use methods and results from \cite{Ho}, \cite{M4}, \cite{AL}, and \cite{SbM}.

Throughout the paper, $G$ is a nontrivial  compact and connected Abelian group
with the normalized Haar measure $m$ and totally ordered character
group $X$, and $X_+$ is the positive cone in $X$. It means that
in the group $X$ there is a distinguished subsemigroup $X_+$ containing the identity
character $1$ and such that $X_+\cap X_+^{-1} = \{1\}$ and $X_+\cup X_+^{-1} = X$. The
rule $\xi\leq\chi:=\xi^{-1}\chi\in X_+$
defines a total order in $X$ which agrees with the structure of the group. We put also $X_-:=X\setminus X_+ (=X_+^{-1}\setminus\{1\})$. As an example,  let $X$ be an
additive subgroup of $\mathbb{R}^n$ endowed with the discrete topology and lexicographical   order so that $G$ is its Bohr compactification (or $X =\mathbb{ Z}^\infty$, the direct sum of countably many copies of $\mathbb{ Z}$ so that $G = \mathbb{T}^\infty$,   the infinite-dimensional torus).
As is well known, a (discrete) Abelian group $X$ can be totally ordered if and
only if it is torsion-free (see, for example, \cite{2}), which in turn is equivalent to the
condition that its character group $G$ is connected (see \cite{3}); the total order on $X$
here is not, in general, unique.

Let $H^p(G)\ (1\leq p\leq\infty)$ denotes the subspace of functions $f\in L^p(G)$ whose
Fourier transforms  $\widehat{f}$ are concentrated on  $X_+$.   We equip  the space   $H^2(G)$ with the inner product $\langle \cdot, \cdot \rangle$ induced from  $L^2(G)$ (see \cite{2}). Note that the set $X$ ($X_+$) is an orthonormal basis of the space $L^2(G)$ (of the space  $H^2(G)$, respectively).

A \textit{Toeplitz operator} $T_\varphi$ on $H^2(G)$ with a symbol
$\varphi\in L^\infty(G)$ is defined as follows:
$T_\varphi f = P_+(\varphi f),\ f\in H^2(G)$ where $P_+:L^2(G)\to H^2(G)$ is the orthogonal  projection.
 These operators   were
defined by Murphy in \cite{M4} and
intensively studied (see, for example, \cite{M4} -- \cite{11}, and \cite{12}).
In the paper \cite{SbM}  Toeplitz operators on Banach spaces $H^p(G)$ were   considered.

\section{$\lambda$-Toeplitz operators over compact Abelian groups and their connection with Toeplitz operators}
\label{2}

\textbf{Definition 2.1.} Let $\lambda:X_+\to\mathbb{C}$. A linear bounded operator $T$ on $H^2(G)$ is said to be  \textit{$\lambda$-Toeplitz} if
$$
\langle T(\chi\xi),\chi\eta\rangle=\lambda(\chi)\langle T\xi,\eta\rangle
$$
for all $\chi,\xi,\eta\in X_+$.

Recall that a \textit{bounded semicharacter} of a discrete semigroup $S$ is a nonzero homomorphism from $S$ into the multiplicative semigroup $\overline{\mathbb{D}}:=\{z\in \mathbb{C}: |z|\leq 1\}$.

\textbf{Lemma 2.2.} \textit{For every $\lambda$-Toeplitz operator  $T\ne 0$ the map $\lambda$ is a bounded semicharacter of the semigroup $X_+$. Moreover, this map extends uniquely to a character of the group $X$ if $|\lambda(\chi)|=1$ for all $\chi\in X_+$. }

\textbf{Proof.}  From the identity
$$
\lambda(\chi_1\chi_2)\langle T\xi,\eta\rangle=\langle T(\chi_1\chi_2\xi),\chi_1\chi_2\eta\rangle=\lambda(\chi_1)\lambda(\chi_2)\langle T\xi,\eta\rangle,
$$
and $\lambda(1)=1$ it follows that $\lambda$ is a nontrivial homomorphism from $X_+$ into the multiplicative semigroup $\mathbb{C}$. Its boundedness follows from the  boundedness  of $T$.
 Now,  let $|\lambda(\chi)|=1$  for all $\chi\in X_+$. Every $\theta\in X$ has the form $\xi^{-1}\chi$ where $\xi, \chi\in X_+$, and   it is easy to verify that $\lambda$ correctly extends  to the character $\lambda$ of the whole group $X$ by the formula $\lambda(\theta):=\lambda(\xi)^{-1}\lambda(\chi)$ . $\Box$

 In the following unless  otherwise stipulated we suppose that  $|\lambda(\chi)|=1$ for all $\chi\in X_+$. Taking  the Pontrjagin-van Kampen duality theorem and  Lemma 2.2 into account we can assume that $\lambda\in G$ and identify $\lambda(\chi)$ with  $\chi(\lambda)$ for $\chi\in X$.

To state our first result we need several definitions.

\textbf{Definition 2.3.} The function $\varphi\in L^2(G)$ with the Fourier transform
$$
\widehat{\varphi}(\chi)=\left\{
\begin{array}
{@{\,}r@{\quad}l@{}}
\langle T1,\chi\rangle  \quad {\rm if }\  \chi\in X_+,\\
\langle T\chi^{-1},1\rangle \quad  {\rm if } \ \chi\in X_-
\end{array}
\right.
$$
is called \textit{the symbol  of $\lambda$-Toeplitz operator} $T$.

By the Bessel's inequality  and Plancherel theorem the above definition is correct because
$$
\sum\limits_{\chi\in X}|\widehat{\varphi}(\chi)|^2=\sum\limits_{\chi\in X_+}|\langle T1,\chi\rangle|^2+\sum\limits_{\xi\in X_+\setminus\{1\}}|\langle T\xi,1\rangle|^2=
$$
$$
\sum\limits_{\chi\in X_+}|\langle T1,\chi\rangle|^2+\sum\limits_{\xi\in X_+\setminus\{1\}}|\langle T^*1,\xi\rangle|^2\leq\|T1\|^2+\|T^*1\|^2.
$$

From now on $T_{\lambda,\varphi}$ denotes the $\lambda$-Toeplitz operator with $\lambda\in G$ and the symbol $\varphi$.

\textbf{Definition 2.4.} By the \textit{modified symbol  of $\lambda$-Toeplitz operator} $T_{\lambda,\varphi}$ we call the function $\varphi_\lambda\in L^2(G)$ with the Fourier transform \footnote{In \cite{Ho} this function was denoted by $\varphi_{\bar{\lambda}, +}$.}
$$
\widehat{\varphi_\lambda}(\chi)=\left\{
\begin{array}
{@{\,}r@{\quad}l@{}}
\overline{\lambda(\chi)}\widehat{\varphi}(\chi) \quad {\rm if }\  \chi\in X_+,\\
\widehat{\varphi}(\chi)  \quad {\rm if } \ \chi\in X_-.
\end{array}
\right.
$$

Define also the unitary operator $U_\lambda$ on  $H^2(G)$  by the formula  $U_\lambda f(x):=f(\lambda x)$.

\textbf{Theorem 2.5.}\textit{The modified symbol $\varphi_\lambda$ of $\lambda$-Toeplitz operator $T_{\lambda,\varphi}$ belongs to $L^\infty(G)$ and}
$$
T_{\lambda,\varphi}=U_\lambda T_{\varphi_\lambda}.
$$
\textit{In particular, $\|T_{\lambda,\varphi}\|=\|\varphi_\lambda\|_\infty$.
}

\textit{Conversely, for every  $\lambda\in G$ and $\psi\in L^\infty(G)$ the operator $U_\lambda T_{\psi}$ is $\lambda$-Toeplitz with modified symbol $\psi$.
}

\textbf{Proof.}  Consider the bounded operator $A=U_\lambda^{-1} T_{\lambda,\varphi}$ on $H^2(G)$. For all $\chi_1,\chi_2\in X_+$ we have
$$
\langle A\chi_1,\chi_2\rangle=\langle U_\lambda^{-1} T_{\lambda,\varphi}\chi_1,\chi_2\rangle=
\langle  T_{\lambda,\varphi}\chi_1,U_\lambda\chi_2\rangle=\overline{\lambda(\chi_2)}\langle  T_{\lambda,\varphi}\chi_1,\chi_2\rangle.
$$

Two cases are possible.

1) $\chi_1^{-1}\chi_2\in X_+$. Then
$$
\langle A\chi_1,\chi_2\rangle=\overline{\lambda(\chi_2)}\lambda(\chi_1^{-1}\chi_2)^{-1}\langle  T_{\lambda,\varphi}(\chi_1^{-1}\chi_2)\chi_1,(\chi_1^{-1}\chi_2)\chi_2\rangle=
$$
$$
\lambda(\chi_1^{-1}\chi_2)^{-1}\overline{\lambda(\chi_2)}\langle  T_{\lambda,\varphi}\chi_2,(\chi_1^{-1}\chi_2)\chi_2\rangle=\overline{\lambda(\chi_1^{-1}\chi_2)}\langle  T_{\lambda,\varphi}1,\chi_1^{-1}\chi_2\rangle.
$$

2) $\chi_1^{-1}\chi_2\in X_-$. Then
$$
\langle A\chi_1,\chi_2\rangle=\overline{\lambda(\chi_2)}\langle  T_{\lambda,\varphi}\chi_2(\chi_1\chi_2^{-1}),\chi_2\rangle=\langle  T_{\lambda,\varphi}(\chi_1\chi_2^{-1}),1\rangle.
$$

So in both cases $\langle A\chi_1,\chi_2\rangle=a(\chi_1^{-1}\chi_2)$ for all $\chi_1,\chi_2\in X_+$   where
$$
a(\chi)=\left\{
\begin{array}
{@{\,}r@{\quad}l@{}}
\overline{\lambda(\chi)}\langle  T_{\lambda,\varphi}1,\chi\rangle \quad {\rm if }\  \chi\in X_+,\\
\langle T_{\lambda,\varphi}\chi^{-1},1\rangle  \quad {\rm if } \ \chi\in X_-.
\end{array}
\right.
$$
 By \cite[Theorem 1]{SbM}, we have  $A=T_\psi$ with $\psi\in L^{\infty}$ and
 $\widehat{\psi}=a$,  $\|T_\psi\|=\|\psi\|_\infty$. Since $a=\widehat{\varphi_\lambda}$, the first statement is proved.

Finally, if   $T=U_\lambda T_{\psi}$ for some $\lambda\in G,  \psi\in L^{\infty}$, then for all $\chi, \xi, \eta\in X_+$ we have, since $U_\lambda^{-1}\zeta=\overline{\lambda(\zeta)}\zeta\ (\zeta\in X_+)$,
$$
\langle T(\chi\xi),\chi\eta\rangle=\langle T_\psi(\chi\xi),U_\lambda^{-1}(\chi\eta)\rangle=\lambda(\chi)\lambda(\eta)\langle T_\psi(\chi\xi),\chi\eta\rangle=
$$
$$\lambda(\chi)\lambda(\eta)\langle T_\psi\xi,\eta\rangle
=\lambda(\chi)\langle T_\psi\xi,U_\lambda^{-1}\eta\rangle=\lambda(\chi)\langle T\xi,\eta\rangle.
$$
Thus $T=T_{\lambda,\varphi}=U_\lambda T_{\varphi_\lambda}$ for some $\varphi_\lambda\in L^\infty(G)$. So $U_\lambda T_\psi=U_\lambda T_{\varphi_\lambda}$, which implies $\psi=\varphi_\lambda$ (see \cite[Theorem 3.2]{M4}).  The proof is complete. $\Box$

In the following,  $S_{\chi}$ ($\chi\in X_+$) denotes the isometry $f\mapsto \chi f$ of the space $H^2(G)$.

\textbf{Corollary 2.6.} \textit{A linear bounded operator $T$ on $H^2(G)$ is   $\lambda$-Toeplitz if and only if}
$$
S_{\chi}^*TS_{\chi}=\lambda(\chi) T\eqno(1)
$$
\textit {for all $\chi\in X_+$.
}

\textbf{Proof.}   The necessity can be verified directly. To prove the sufficiency, consider the operator $A:=U_\lambda^{-1} T$. If we replace $T$ by $U_\lambda A$ in (1), we get
for all $\chi, \xi, \eta\in X_+$
$$
\langle S_{\chi}^*U_\lambda A S_{\chi}\xi, \eta\rangle=\lambda(\chi)\langle U_\lambda A \xi, \eta\rangle.
$$
In other words,
$$
\langle  A(\chi\xi), U_\lambda^{-1}(\chi\eta)\rangle=\lambda(\chi)\langle A\xi, U_\lambda^{-1}\eta\rangle.
$$
Using $U_\lambda^{-1}\zeta=\overline{\lambda(\zeta)}\zeta\ (\zeta\in X_+)$, we obtain
$$
\langle  A(\chi\xi), \chi\eta\rangle=\langle A\xi,\eta\rangle.
$$
It  follows \cite[Theorem 3.10]{M4} that $A=T_\psi$ for some $\psi\in L^\infty(G).\ \Box$

\textbf{Remark.}  As it was mentioned in \cite{Ho}  for the classical case $G=\mathbb{T}$ the problem of studing  operators with the property (1)  was posed in \cite[p. 629 -- 630]{BH}.

\textbf{Corollary 2.7.} \textit{Every $\lambda$-Toeplitz operator is uniquely determined by the pair $(\lambda, \varphi)$ where $\lambda\in G$ and $\varphi \in L^2(G)$ is such that $\varphi_\lambda\in L^\infty(G)$.}

\textbf{Proof.}  The necessity for the condition $\varphi_\lambda\in L^\infty(G)$ was proved above. Conversely, if the pair  $(\lambda, \varphi)$ meets all the conditions of this corollary, the operator $T:=U_\lambda T_{\varphi_\lambda}$ is $\lambda$-Toeplitz with a symbol $\varphi$ and it is obvious in view of the preceding theorem that two $\lambda$-Toeplitz operators with the same  symbol  $\varphi$ coincide. $\Box$

Let $C(G)$ denotes the algebra  of continuous functions on $G$ and $C(G)^{-1}$  the group  of  invertible elements of  $C(G)$. To state and prove the next corollary (and several other results below), we need  the notion of  the rotation index for functions in some subgroup of
$C(G)^{-1}$  given in \cite{SbM}. We begin with the definition
of the rotation index for a character of the group $G$ (the symbol $\sharp F$ will denote the
number of elements of a finite set $F$ in what follows).

\textbf{Definition 2.8.} In each of the following cases, we define \textit{the rotation index of a character}
$\chi\in X$ as follows:

1) ${\rm ind}\chi= \sharp(X_+ \setminus\chi X_+)$ if $\chi\in  X_+$ and the set $X_+ \setminus\chi X_+$ is finite;

2) ${\rm ind}\chi= {\rm ind}\chi_1-{\rm ind}\chi_2$ if $\chi=\chi_1\chi_2^{-1}$, where $\chi_j\in X_+$, where both sets $X_+ \setminus\chi_j X_+$
are finite, $j = 1, 2$.

In the other cases we assume that the character has no index.

We denote the set of characters that have an index by $X^i$. It follows from  \cite[Theorem 2]{SbM}  that $X^i$ is a cyclic subgroup of $X$ and it is nontrivial if and only if $X$ contains the smallest strictly positive element.

\textbf{Definition 2.9.} Consider a function $\varphi\in C(G)^{-1}$ of the form $\chi \exp(g)$, where $g \in C(G)$ and $\chi\in  X$  (the Bohr-van Kampen decomposition). If $\chi\in  X^i$, then we set
${\rm ind}\varphi= {\rm ind}\chi$.
Otherwise we assume that the function $\varphi$ has no index.

We denote the set of functions in $C(G)^{-1}$ which have an index by $\Phi(G)$. Thus, $\Phi(G)=X^i\exp(C(G))$.

We recall that, for a bounded operator $T$ on a Banach space $Y$,
the symbol $T \in \Phi_+(Y)$  means that the image ${\rm Im} T$ is closed and the kernel ${\rm Ker T}$ is
finite-dimensional, whereas the symbol  $T \in \Phi_-(Y)$  means that the quotient space
$Y/ {\rm Im} T$ is finite-dimensional; the operators in  $\Phi_+(Y)\cup \Phi_-(Y)$ are referred to as
semi-Fredholm operators, whereas those in $\Phi_+(Y)\cap \Phi_-(Y)$ are called Fredholm
operators on $Y$.

\textbf{Corollary 2.10.} \textit{Let $\varphi_\lambda\in  C(G)$. The $\lambda$-Toeplitz operator $T_{\lambda,\varphi}$ is Fredholm if and only if $\varphi_\lambda\in  \Phi(G)$. In this case,
$$
{\rm Ind}T_{\lambda,\varphi}=-{\rm ind} \varphi_\lambda.
$$
}
\textbf{Proof.} It follows from the above theorem and \cite[Theorem 4]{SbM}. $\Box$

\textbf{Corollary 2.11.}  \textit{If an operator $T_{\lambda,\varphi}$ is semi-Fredholm, then its
modified symbol $\varphi_\lambda$ is invertible in the algebra $L^\infty(G)$.}

\textbf{Proof.} It follows from the above theorem and  \cite[Theorem 3]{SbM}. $\Box$

\textbf{Corollary 2.12.} \textit{A $\lambda$-Toeplitz operator is compact if and only if it is zero.}

\textbf{Proof.} It follows from the above theorem and  \cite[Theorem 3.5]{M4}. $\Box$

\textbf{Corollary 2.13.} \textit{A linear bounded operator $T$ on $H^2(G)$ is   $\lambda$-Toeplitz if and only if} $\overline{\chi_2(\lambda)}\langle T\chi_1,\chi_2\rangle$ \textit {depends on $\chi_1^{-1}\chi_2$ only} ($\chi_1, \chi_2\in X_+$).

\textbf{Proof.} It follows from the above theorem and  \cite[Theorem 1]{SbM}. $\Box$

\section{Rotational invariance for the essential
spectrum}
\label{3}

In the following,  $\sigma_e(T)$ stands for the \textit{essential (Fredholm) spectrum } of an operator $T$,  $\rho_e(T)=\mathbb{C}\setminus \sigma_e(T)$.

\textbf{Lemma 3.1.} \textit{Let $\lambda$ be a nonvanishing complex-valued function defined on the set $X_+\cap X^i$, $T$ a linear bounded operator on $H^2(G)$ such that $S_{\chi}^*TS_{\chi}=\lambda(\chi) T$
 for all $\chi\in X_+\cap X^i$. Suppose that $T-\mu$ is Fredholm. Then for every  $\chi\in X_+\cap X^i$ the operator $T-\lambda(\chi)^{-1}\mu$ is also Fredholm  and}
$$
{\rm Ind}(T-\mu)={\rm Ind}(T-\lambda(\chi)^{-1}\mu).
$$
\textit{Consequently, for every  $\chi\in X_+\cap X^i$ we have $\lambda(\chi)^{-1}\rho_e(T)\subseteq \rho_e(T)$.}

\textbf{Proof.} First note that $S_\chi (=T_\chi)$ is Fredholm for every $\chi\in X_+\cap X^i$ by \cite[Theorem 4]{SbM}.
Now the operator $T-\lambda(\chi)^{-1}\mu$ is  Fredholm, since
$$
T-\lambda(\chi)^{-1}\mu=\lambda(\chi)^{-1}S^*_\chi(T-\mu)S_\chi.
$$
We conclude  also that $\rho_e(T)\subseteq \lambda(\chi)\rho_e(T)$, and
$$
{\rm Ind}(T-\lambda(\chi)^{-1}\mu)={\rm Ind}S^*_\chi+{\rm Ind}(T-\mu)+{\rm Ind}S_\chi={\rm Ind}(T-\mu). \Box
$$

By definition, put  $\tau(x) = \lambda^{-1}x,\ x\in G$.

\textbf{Lemma 3.2.} \textit{ For any $\lambda$-Toeplitz operator $T_{\lambda,\varphi}$ we have
$$
T_{\lambda,\varphi}^k=U_\lambda^kT_{\varphi_\lambda\circ\tau^{k-1}}\cdots T_{\varphi_\lambda}
$$
for $k=1,2,\dots$.
}

\textbf{Proof.}  First we prove  that
$$
T_fU_\lambda^k=U_\lambda^k T_{f\circ \tau^k}
$$
for $f\in L^\infty(G), k=1,2,\dots$.

Indeed, using \cite[Theorem 1]{SbM} we have for $\chi, \eta \in X_+$
$$
\langle T_fU_\lambda^k \chi,\eta\rangle=\lambda(\chi)^k\langle T_f\chi,\eta\rangle=\lambda(\chi)^k\widehat{f}(\chi^{-1}\eta),
$$
and, on the other hand,
$$
\langle U_\lambda^k T_{f\circ \tau^k}\chi,\eta\rangle=\langle T_{f\circ \tau^k}\chi, U_{\lambda^{-k}}\eta\rangle=\lambda(\eta)^k\langle T_{f\circ \tau^k}\chi, \eta\rangle=
$$
$$
\lambda(\eta)^k\widehat{f\circ \tau^k}(\chi^{-1}\eta)=\lambda(\chi)^k\widehat{f}(\chi^{-1}\eta).
$$
Now in view of Theorem 2.5 the statement of the Lemma follows by induction. $\Box$

In what follows we put\footnote{In \cite{Ho} this function was denoted by $\Phi_{\bar{\lambda}+}$.}  $\Phi_\lambda=\prod_{j=0}^{q-1}\varphi_\lambda\circ\tau^j$.

Our main result is the following.

\textbf{Theorem 3.3.} \textit{Let $\lambda\in G$.}

1)  \textit{$\sigma_e(T_{\lambda,\varphi})=\lambda(\chi)\sigma_e(T_{\lambda,\varphi})$ for all $\chi\in X^i$.}

2)\textit{ If the number $\lambda(\chi)$ is not of finite order in  $\mathbb{T}$ for some $\chi\in X^i$, then  $\sigma_e(T_{\lambda,\varphi})$ is circular.}

3)  \textit{Suppose  $\lambda$ is of order $q$ in $G$, and the number $\lambda(\chi_0)$ is a primitive root of $1$ of order $q$ for some $\chi_0\in X^i_+$. Then}

${\rm Ind}(T_{\lambda,\varphi}-\mu)=q^{-1}{\rm Ind}(T_{\varphi_\lambda\circ\tau^{q-1}}\dots T_{\varphi_\lambda}-\mu^q)$ \textit{for all} $\mu\in\rho_e(T_{\lambda,\varphi})$.

 4) \textit{Let} $\varphi\in H^\infty(G)$. \textit{Then $T_{\lambda,\varphi}$  is a weighted shift operator on $H^2(G)$ and}

 $\sigma(T_{\lambda,\varphi})=\lambda(\chi)\sigma(T_{\lambda,\varphi})$ \textit{for all} $\chi\in X$.
 \textit{Consequently,  $\sigma(T_{\lambda,\varphi})$ is circular if $\lambda$ is not of finite order in  $G$.}

5) \textit{Let} $\varphi\in H^\infty(G)$. \textit{Suppose  $\lambda$ is of order $q$ in $G$ and the number $\lambda(\chi_0)$ is a primitive root of $1$ of order $q$ for some $\chi_0\in X^i$. Then
}

(i) $\sigma_e(T_{\lambda,\varphi})=\{\mu\in \mathbb{C}:\mu^q\in \sigma_e(T_{\Phi_\lambda})\}$;

(ii) $\sigma(T_{\lambda,\varphi})=\{\mu\in \mathbb{C}:\mu^q\in \sigma(T_{\Phi_\lambda})\}.$

\textbf{Proof.} To prove 1),  we can assume that $X^i\ne \{1\}$. First suppose that the number $\lambda(\chi)$ is  of finite order  in  $\mathbb{T}$  for some $\chi\in X^i\setminus\{1\}$. It follows, since the group $X^i$ is cyclic \cite[Theorem 2]{SbM}, that the number  $\lambda(\chi)$ is  of finite order   for every $\chi\in X^i$. Let $\chi\in X_+\cap X^i$ and  $\lambda(\chi)$ is  of  order $r$. Then  Corollary 2.6 and Lemma 3.1 imply
$$
\rho_e(T_{\lambda,\varphi})\supseteq \lambda(\chi)^{-1}\rho_e(T_{\lambda,\varphi})\supseteq\dots\supseteq \lambda(\chi)^{-r}\rho_e(T_{\lambda,\varphi})=\rho_e(T_{\lambda,\varphi}).
$$
Since, by  \cite[Theorem 2, 2)]{SbM},  $X^i=( X_+\cap X^i)\cup ( X_+\cap X^i)^{-1}$, this proves 1) in the case of finite order.

Now suppose that the number $\lambda(\chi)$ is  not of finite order  in  $\mathbb{T}$  for some (and therefore for all) $\chi\in X^i\setminus \{1\}$ and choose $\chi_1\in  ( X_+\cap X^i)\setminus \{1\}$. Let $\mu\in \rho_e(T_{\lambda,\varphi})$. There is an arc $J$ in the circle $\{z:|z|=|\mu|\}$ such  that $\mu\in J\subset \rho_e(T_{\lambda,\varphi})$. Since $\{\overline{\lambda(\chi_1)}^k:k=1,2\dots\}$ is dense in the circle $\mathbb{T}=\{z:|z|=1\}$, Lemma 3.1 yields
$$
 \rho_e(T_{\lambda,\varphi})\supseteq \bigcup\limits_{k=1}^\infty\overline{\lambda(\chi_1)}^kJ=\{z:|z|=|\mu|\}.\eqno(2)
$$
This proves 2). Moreover, formula (2) implies  $\lambda(\chi)\mu\in \rho_e(T_{\lambda,\varphi})$  and therefore   $\rho_e(T_{\lambda,\varphi})\subseteq \overline{\lambda(\chi)}\rho_e(T_{\lambda,\varphi})$ for all $\chi\in X^i$. Combining this result and Lemma 3.1, we obtain the assertion 1) for $\chi\in X_+\cap X^i$ and therefore for all $\chi\in X^i$. Thus, 1) is proved in the case of infinite order as well.

3) Consider the factorization
$$
T_{\lambda,\varphi}^q-\mu^q=\prod\limits_{k=0}^{q-1}(T_{\lambda,\varphi}-\overline{\lambda(\chi_0)}^k\mu).
$$
  By 1), $\overline{\lambda(\chi_0)}^k\mu\in \rho_e(T_{\lambda,\varphi})$.  Hence   the operator $T_{\lambda,\varphi}^q-\mu^q$ is Fredholm. Since, by Lemma 3.2, $T_{\lambda,\varphi}^q=T_{\varphi_\lambda\circ\tau^{q-1}}\dots T_{\varphi_\lambda}$, we conclude in view of Lemma 3.1 and Corollary 2.6 that
  $$
  {\rm Ind}(T_{\varphi_\lambda\circ\tau^{q-1}}\dots T_{\varphi_\lambda}-\mu^q)=\sum\limits_{k=0}^{q-1}{\rm Ind}(T_{\lambda,\varphi}-\overline{\lambda(\chi_0^k)}\mu)=q{\rm Ind}(T_{\lambda,\varphi}-\mu).
  $$

4) First note that, by Theorem 2.5, for $\varphi\in H^\infty(G)$ the $\lambda$-Toeplitz operator has the form
$$
T_{\lambda,\varphi}f(x)=\varphi(x)f(\lambda x)\eqno(3)
$$
and so is a weighted shift operator on $H^2(G)$ (we refer the reader to \cite{AL} for the general theory of weighted shift operators).  It follows that for all $\chi\in X$
$$
S^{-1}_\chi T_{\lambda,\varphi}S_\chi=\chi(\lambda)T_{\lambda,\varphi}
$$
(actually, both sides of this  equality are bounded operators and coincide on the basis $X_+$ of $H^2(G)$). In turn, the last equality implies
$$
\overline{\chi(\lambda)}\mu-T_{\lambda,\varphi}=\overline{\chi(\lambda)}S^{-1}_\chi (\mu-T_{\lambda,\varphi})S_\chi.
$$
This proves the first statement. Now the  second statement follows from the fact that the set $\{\lambda(\chi):\chi\in X\}$ is dense in $\mathbb{T}$ if $\lambda$ is an element of infinite order  (see, e. g., \cite[p. 119]{AL}).

5)
(i) Since $\varphi \in H^\infty(G)$, we have  $\varphi_\lambda\circ \tau^j \in H^\infty(G)$ and therefore $T_{\lambda,\varphi}^q=T_{\Phi_\lambda}$ by Lemma 3.2. Hence, by the spectral mapping theorem,
$$
\sigma_e(T_{\lambda,\varphi})^q=\sigma_e(T_{\Phi_\lambda}), \mbox{ and } \sigma(T_{\lambda,\varphi})^q=\sigma(T_{\Phi_\lambda}).\eqno(4)
$$
 First of  all, (4) implies that
$$
 \sigma_e(T_{\lambda,\varphi})\subseteq\{\mu\in \mathbb{C}:\mu^q\in \sigma_e(T_{\Phi_\lambda})\}.
$$
Now let $\mu\in \mathbb{C}, \mu^q\in \sigma_e(T_{\Phi_\lambda})$. Then, by (4), one can find $\nu\in \sigma_e(T_{\Phi_\lambda})$ such that $\nu^q=\mu^q$. Using 1),  we get  for some $k\in\{1,\dots,q\}$ that $\mu=\nu\lambda(\chi_0)^k\in \sigma_e(T_{\lambda,\varphi})$. This proves (i).

(ii) In view of 4) the  proof of  this equality  is similar to that of (i). $\Box$

\textbf{Remark.} For $f\in C(G)$ spectra $\sigma_e(T_f)$ and $\sigma(T_f)$ were calculated in \cite{SbM}. For  $f\in H^\infty(G)$  the spectrum $\sigma(T_f)$ equals to $\sigma_{H^\infty(G)}(f)$, the spectrum of the element $f$ in the algebra $H^\infty(G)$ \cite[Theorem 3.12]{M4}.

\textbf{Remark.} Without the assumption that $\lambda(\chi_0)$ is a primitive root of $1$ of order $q$ the conclusions of the part  5) of the above theorem are false as the following simple example shows.
 Let $G=\mathbb{T}, \lambda =1, q=2, \varphi(z)=z+2$. Then $T_{\lambda,\varphi}=T_{\varphi}$ is an analytic Toeplitz operator and therefore   $\sigma_e(T_{\lambda,\varphi})=\varphi(\mathbb{T})=\mathbb{T}+2,  \sigma(T_{\lambda,\varphi})={\rm cl}(\varphi(\mathbb{D}))=\overline{\mathbb{D}}+2$ (see, e. g.,  \cite{Pel}, p. 98, p. 93). On the other hand, here $\Phi_\lambda(z)=(z+2)^2$. It follows that
$\sigma_e(T_{\lambda,\varphi})\ne\{\mu\in \mathbb{C}:\mu^2\in \Phi_\lambda(\mathbb{T})\}$ and
$\sigma(T_{\lambda,\varphi})\ne\{\mu\in \mathbb{C}:\mu^2\in {\rm cl}(\Phi_\lambda(\mathbb{D}))\}$.

\section{$\lambda$-Toeplitz operators with Arens - Singer symbols}
\label{4}

\textbf{Definition 4.1.} The uniform algebra $A(G):=H^\infty(G)\cap  C(G)$ is called \textit{the Arens - Singer algebra } (of the group  $G$).

We have the following corollary of Theorem 3.3.

\textbf{Corollary 4.2.} \textit{Let $\varphi_\lambda\in A(G)$. Suppose  $\lambda$ is of order $q$ in $G$, and the number $\lambda(\chi_0)$ is a primitive root of $1$ of order $q$ for some $\chi_0\in X^i_+$. Then}
 $$
 {\rm Ind}(T_{\lambda,\varphi}-\mu)=-q^{-1}{\rm ind}(\Phi_\lambda-\mu^q)
 $$
 \textit{for all} $\mu\in\rho_e(T_{\lambda,\varphi})$;

\textbf{Proof.}  It follows from the statement 3) of Theorem 3.3 and \cite[Theorem 4]{SbM}, since $T_{\varphi_\lambda\circ\tau^{q-1}}\dots T_{\varphi_\lambda}-\mu^q=T_{\Phi_\lambda-\mu^q}$ for  $\varphi_\lambda\in A(G)$.  $\Box$

Recall that a topological group $G$ is called \textit{monothetic with a generator} $\lambda$ if the set $\{\lambda^n:n\in \mathbb{Z}\}$ is dense in $G$. For  example, the  tori  $\mathbb{T}^\mathfrak{m}, \mathfrak{m}\leq \aleph_0$ are  (compact and connected) monothetic groups with a countable base of the topology.

The following  fact is a corollary of \cite[Theorem 4.4, and Theorem 5.22]{AL}.

\textbf{Lemma 4.3.} \textit{Let $G$ be monothetic with a countable base of the topology and $\lambda$ a generator of $G$. If $\varphi\in A(G)$, then $r(T_{\lambda,\varphi})$, the spectral radius of $T_{\lambda,\varphi}$, is}
$$
\exp\int\limits_G\log|\varphi|dm.
$$
\textbf{Proof.} Since $T_{\lambda,\varphi}$  is a weighted shift operator, we have by  \cite[Theorem 4.4]{AL}
$$
r(T_{\lambda,\varphi})=\max\limits_{\nu\in \Lambda}\exp\int\limits_G\log|\varphi|d\nu,
$$
where $\Lambda$ denotes the set of all  measures on the Shilov boundary $\partial A(G)$ of the algebra $A=A(G)$, which are ergodic with respect to some homeomorphism $\alpha$ of $\partial A(G)$ associated with $T_{\lambda,\varphi}$. But it is known that $\partial A(G)=G$ \cite{AS} and it is easy to verify that $\alpha(x)=\lambda x$ for all $x\in G$. By \cite[Theorem 5.22]{AL}, for $\nu\in \Lambda$ there is such $x_0\in G$ that
$$
\int\limits_Gfd\nu=\int\limits_Gf(gx_0)dm(g)=\int\limits_Gfdm
$$
for all $f\in C(G)$. Thus, $\Lambda=\{m\}$ which completes the proof. $\Box$

The next theorem is a partial generalization of   \cite[Proposition 3.5]{Ho}.

\textbf{Theorem 4.4.}  \textit{Let $G$ be monothetic with a countable base of the topology and $\lambda$ a generator of $G$. Suppose that $\varphi\in A(G)\cap C(G)^{-1}$}.

1)  \textit{If  $\varphi$ is invertible in $A(G)$, then}
$$
\sigma(T_{\lambda,\varphi})=\sigma_e(T_{\lambda,\varphi})=\{\mu:|\mu|=r(T_{\lambda,\varphi})\}.
$$

2)  \textit{If $\varphi$ is not invertible in $A(G)$,  then}
$$
\sigma(T_{\lambda,\varphi})=\{\mu:|\mu|\leq r(T_{\lambda,\varphi})\}.
$$

\textbf{Proof.} 1) Let $\varphi^{-1}\in A(G)$. Since  $\varphi$ is  invertible in $H^\infty(G)$,
the analytic Toeplitz operator $T_\varphi$ is  invertible and   $T_\varphi^{-1}=T_{\varphi^{-1}}$.
Formula (3) implies that $T_{\lambda,\varphi}=T_{\varphi}U_\lambda$. Whence, this operator is invertible and
$$
T_{\lambda,\varphi}^{-1}=U_{\lambda^{-1}}T_{\varphi^{-1}}=T_{\lambda^{-1},\varphi^{-1}\circ\tau}.
$$
 Moreover, by Lemma 4.3,
$$
r(T_{\lambda,\varphi}^{-1})=\exp\int\limits_G\log|\varphi^{-1}(\lambda^{-1} x)|dm(x)=\frac{1}{r(T_{\lambda,\varphi})}.
$$
For  $0<|\mu|<r(T_{\lambda,\varphi})$ we have
$$
T_{\lambda,\varphi}-\mu=\mu T_{\lambda,\varphi}(\mu^{-1}-T_{\lambda,\varphi}^{-1}).
$$
It follows that the operator $T_{\lambda,\varphi}-\mu$ is  invertible, since $|\mu^{-1}|>1/r(T_{\lambda,\varphi})=r(T_{\lambda,\varphi}^{-1})$. But the spectrum of  $T_{\lambda,\varphi}$ is circular by the assertion 4) of Theorem 3.3, since $\lambda$ is not of finite order in $G$. This proves that $\sigma(T_{\lambda,\varphi})$ is a circle centered at the origin.

Now let us suppose that $\sigma(T_{\lambda,\varphi})\ne \sigma_e(T_{\lambda,\varphi})$.  Then the complement in $\mathbb{C}$ of $\sigma_e(T_{\lambda,\varphi})$  is connected, and  by \cite[Corollary XI.8.5 and Theorem II.1.1]{GGK} the points of  $\sigma(T_{\lambda,\varphi})\setminus \sigma_e(T_{\lambda,\varphi})$ are isolated points of $\sigma(T_{\lambda,\varphi})$. This contradiction concludes the  proof of the first statement.

2) First let us prove that $0\in \sigma(T_{\lambda,\varphi})$. Indeed, the function  $\varphi$ is not invertible in $H^\infty(G)$, since $\varphi\in C(G)^{-1}$ and $A(G)=H^\infty(G)\cap  C(G)$. It follows that $T_{\varphi}$ is not invertible. Assume the converse. Then $H^2(G)=T_{\varphi}H^2(G)=\varphi H^2(G)$, and therefore $\psi\varphi=1$ for some $\psi\in H^2(G)$. This implies that $H^2(G)=\psi\varphi  H^2(G)=\psi H^2(G)$ and hence  $\psi\in H^\infty(G)$ by \cite[Lemma 2]{SbM}. Thus we arrive at a contradiction. Because of the equality $T_{\lambda,\varphi}=T_{\varphi}U_{\lambda}$,  the operator $T_{\lambda,\varphi}$ is not invertible along with $T_{\varphi}$.

Since, by Theorem 3.3,  the set $\sigma(T_{\lambda,\varphi})$ is circular, it remains to prove that the set $|\sigma(T_{\lambda,\varphi})|:=\{|\mu|:\mu\in \sigma(T_{\lambda,\varphi})\}$ is connected. But it follows from \cite[Corollary 7.4]{AL} that the number of connected components of  $|\sigma(T_{\lambda,\varphi})|$ do not exceed the number of connected components of  $M$, the maximal ideal space  of $A(G)$.  In turn, the space $M$ can be identified with $\widehat{X_+}$, the space of bounded semicharacters of the semigroup $X_+$ \cite{AS}. Since $X_+\cap X_+^{-1}=\{1\}$, the connectedness of $\widehat{X_+}$ follows from \cite[Lemma 3]{MR}. $\Box$

\textbf{Remark.} In the case of $G=\mathbb{T}$ it is known  that under the conditions of the part 2) of Theorem 4.4 $\sigma_e(T_{\lambda,\varphi})$ is a circle centered at the origin, too  \cite[Proposition 3.4]{Ho}. This is not the case for groups distinct from $\mathbb{T}$ because $0\in  \sigma_e(T_{\lambda,\chi})$ if $\chi\in X_+\setminus X^i$ (by \cite[Corollary 1]{SbM}, $X_+\setminus X^i\ne\emptyset$ if $G\ne\mathbb{T}$). In fact, in this case, the operator $T_\chi$ is not Fredholm by \cite[Theorem 4]{SbM}. Hence, the operator $T_{\lambda,\chi}= T_{\chi}U_\lambda$ is not Fredholm as well.




\bibliographystyle{elsarticle-num}


\vspace{5mm}

\end{document}